\documentclass[10pt]{article}

\usepackage{amsmath}
\usepackage{amsopn}
\usepackage{amsthm}
\usepackage[T1]{fontenc}
\usepackage{amssymb}
\usepackage[frenchb,english]{babel}
\usepackage[mathscr]{eucal}
\usepackage{amsbsy}

\theoremstyle{plain}
\newtheorem{theoreme}{Th\'{e}or\`{e}me}[section]
\newtheorem{definition}[theoreme]{D\'{e}finition}

\newtheorem{corollaire}[theoreme]{Corollaire}

\def\one{\hbox{1\hskip -3pt I}}
\def\rd{\Bbb{R}^d}
\def\sd{S^{d-1}}
\def\qed{\hspace*{2mm} \hfill $\Box$\bigskip}
\def\proof{\par\rm{\bf D\'emonstration:}\quad}

\def\reip{{\bf P}}

\def\reie{{\bf E}}

\def\cvl1{\stackrel{\mathcal{L}^1}{\longrightarrow}}

\title{Transformations des lois multivari\'ees avec queues r\'eguli\`eres}
\author{%
\renewcommand{\thefootnote}{\alph{footnote}}
Youri \textsc{Davydov}\,\footnotemark[1]{}\; , %
Shuyan \textsc{Liu}\,\footnotemark[1]{}
}
\date{}
\begin{document}
\addtocounter{page}{-1}
\maketitle
\renewcommand{\thefootnote}{\alph{footnote}\,}
\footnotetext[1]{Laboratoire P. Painlev\'e,
UMR 8524 CNRS Universit\'e Lille I, B\^at M2, Cit\'e Scientifique,
F-59655 Villeneuve d'Ascq Cedex, France.}

\renewcommand{\thefootnote}{\arabic{footnote}}

\begin{abstract}
Let $X$ be a random vector in $\rd$ with a regularly varying tail. We consider two transformations $\|X\|f(\frac{X}{\|X\|})$, $f: \sd\rightarrow\sd$, and $Xf(\frac{X}{\|X\|})$, $f: \sd\rightarrow \mathbb{R}_+$. Some sufficient conditions for preserving the property of regularity of the tail for this kind of transformations are given.
\end{abstract}

\selectlanguage{frenchb}
\begin{abstract}
Soit $X$ un vecteur al\'eatoire  dans $\rd$ \`a queue \`a variation r\'eguli\`ere. On consid\`ere deux transformations $\|X\|f(\frac{X}{\|X\|})$, $f: \sd\rightarrow\sd$, et $Xf(\frac{X}{\|X\|})$, $f: \sd\rightarrow \mathbb{R}_+$. Nous donnons des conditions suffisantes pour que la propri\'et\'e de r\'egularit\'e de la queue soit pr\'eserv\'ee sous les transformations de ce type.
\end{abstract}



\newpage


\section{Introduction}

Variation r\'eguli\`ere est un des concepts de base qui appara\^it de fa\c con naturelle dans les diff\'erents contextes de la
th\'eorie des probabilit\'es et ses applications. On renvoie le lecteur aux monographies de Feller \cite{feller}, de Araujo \& gin\'e\cite{arauj}, de Resnick \cite{res} et de Bingham et al. \cite{bing} pour la pr\'esentation exhaustive de la mati\`ere. Pour souligner l'importance de cette notion rappelons qu'elle est li\'ee \'etroitement avec la caract\'erisation des domaines d'attraction des lois stables multidimensionnelles (voir Araujo et al. \cite{arauj}, Samorodnitsky et al. \cite{samo} et Davydov et al. \cite{strict} ). On trouve plusieurs informations sur les propri\'et\'es et les applications dans les articles r\'ecents de Mikosch \cite{mikosch2003mda}, Basrak et al. \cite{basrak2002rvg} et Jacobsen et al. \cite{jacobsen}. Rappelons d'abord la d\'efinition.

\begin{definition}\label{def1}
La loi du vecteur al\'eatoire (v.a.) $X$ dans $\rd$ est dite {\em \`a queue \`a variation r\'eguli\`ere} si il existe une mesure finie $\sigma$ sur la sph\`ere unit\'e $\sd$, un nombre $\alpha>0$ et une fonction $L$ \`a variation lente tels que
\begin{equation}\label{regulier}
\lim_{x\rightarrow\infty}\frac{x^\alpha}{L(x)}\reip\left\{\frac{X}{\|X\|}\in B, \|X\|>x\right\}=\sigma(B)
\end{equation}
pour tous $B\in\mathcal{B}(\sd)$ avec $\sigma(\partial B)=0$; ici $\|\cdot\|$ est la norme euclidienne.
\end{definition}
On appelle $\sigma$ {\em mesure spectrale} de $X$, et $\alpha$ s'appelle {\em exposant de variation r\'eguli\`ere}.

Le fait que $X$ est \`a queue \`a variation r\'eguli\`ere sera dans la suite not\'e par l'\'ecriture "$X\in \mbox{VR}(\alpha,\sigma)$". 

Il est clair que sans perdre la g\'en\'eralit\'e on peut consid\'erer $\sigma$ comme une mesure normalis\'ee: $\sigma(\sd)=1$. En prenant dans (\ref{regulier}) $B=\sd$, on d\'eduit imm\'ediatement que 
\begin{equation}\label{regulier1}
\frac{x^\alpha}{L(x)}\reip\{ \|X\|>x\}\rightarrow 1, \  x\rightarrow\infty,
\end{equation}
c'est-\`a-dire que $ \|X\|$ est la variable al\'eatoire positive \`a queue \`a variation r\'eguli\`ere.

Les relations (\ref{regulier}) et (\ref{regulier1}) donnent
\[\reip\left\{\left. \frac{X}{\|X\|}\in B \  \right|  \  \|X\|>x\right\}\rightarrow \sigma(B), \  x\rightarrow\infty,\]
pour tous $B\in\mathcal{B}(\sd)$ avec $\sigma(\partial B)=0$, ce qui signifie que la loi conditionnelle de $\frac{X}{\|X\|}$ sachant $\{\|X\|>x\}$ converge faiblement vers $\sigma$.

Il existent des diff\'erentes caract\'erisations de la propri\'et\'e $X\in \mbox{VR}(\alpha,\sigma)$ (voir, i.e. Mikosch \cite{mikosch2003mda}), on n'en donne ici que deux.

\begin{enumerate}
\item Le v.a. $X\in \mbox{VR}(\alpha,\sigma)$ si et seulement si il existe une fonction $\tilde L$ \`a variation lente telle que pour tous $r>0$ et $B\in\mathcal{B}(\sd)$ avec $\sigma(\partial B)=0$
\begin{equation}\label{regulier2}
\lim_{n\rightarrow\infty}n\reip\left\{\frac{X}{\|X\|}\in B, \|X\|>rb_n\right\}=\sigma(B)r^{-\alpha}
\end{equation}
avec $b_{n}=n^{1/\alpha}\tilde L(n)$.

\item Pour formuler le deuxi\`eme crit\`ere on passe aux coordonn\'ees polaires et identifie $\rd\backslash\{{\bf 0}\}$ avec le produit $(0,\infty)\times\sd$. On introduit les mesures $Q_{n}$ et $Q$ sur $\mathcal{B}((0,\infty)\times\sd)$ par 
\begin{equation}\label{defQn}
Q_n((r,\infty)\times B)=n\reip\left\{\frac{X}{\|X\|}\in B, \|X\|>rb_n\right\},
\end{equation}
\begin{equation}\label{defQ}
Q=m_\alpha\times\sigma,
\end{equation}
o\`u $m_\alpha(dr)=\alpha r^{-\alpha-1}dr$.
Alors $X\in \mbox{VR}(\alpha,\sigma)$ est \'equivalent \`a la convergence vague
\begin{equation}\label{def3}
Q_n \stackrel{vag}\rightarrow Q, \; \; n\rightarrow\infty.
\end{equation}

\end{enumerate}

Soit $\mathcal{S}_{+}$ la famille des fonctions d\'efinies sur $(0,\infty)\times\sd$ dont les supports sont s\'epar\'es de z\'ero, c'est-\`a-dire $f\in \mathcal{S}_{+}$ si et seulement si $\exists \varepsilon >0$ tel que $\mbox{supp}(f)\subset [\varepsilon,\infty)\times\sd$.

Il n'est pas difficile \`a remarquer que la convergence en (\ref{def3}) entraine la convergence
\begin{equation}\label{def4}
\int fdQ_{n}\rightarrow\int fdQ
\end{equation}
pour toute $f\in\mathcal{S}_{+}$ qui est born\'ee et $Q$-p.p. continue.

D'autre part pour \'etablir (\ref{def3}) il suffit de v\'erifier (\ref{def4}) pour toutes $f\in\mathcal{S}_{+}$ continues born\'ees.

La variation r\'eguli\`ere a la propri\'et\'e importante qu'elle est pr\'eserv\'ee sous plusieurs opr\'erations et transformations qu'on utilise souvent en pratique. Une large collection des r\'esultats de ce type est pr\'esent\'ee dans le survey de Jessen et Mikosch \cite{jessen2006rvf}.

Le but de notre travail est de compl\'eter les investigations dans cette direction. On consid\`ere ici deux genres des transformations. Premi\`erement, on \'etudie le passage du v.a. initial $X$ au vecteur $Y=\|X\|f(\frac{X}{\|X\|})$ o\`u $f$ est une application de $\sd$ \`a $\sd$. Dans le deuxi\`eme cas on transforme la partie radiale de $X$, plus exactement on s'int\'eresse \`a la variation r\'eguli\`ere du vecteur $Y=Xf(\frac{X}{\|X\|})$, o\`u la fonction $f$ cette fois-ci est une application de $\sd$ \`a $\mathbb{R}_+$. On propose des conditions suffisantes et on donne des exemples qui montrent que ces conditions ne pourront \^etre affaiblies sensiblement. En conclusion remarquons que les propri\'et\'es des transformations pr\'esent\'ees ici seront utiles pour les simulations des vecteurs appartenant au domaine d'attraction d'une loi stable avec la mesure spectrale donn\'ee. Dans ce contexte on peut mentionner les articles de Chambers et al. \cite{chambers} et de Modarres et al. \cite{modarres}.


\section{R\'esultats}

{\bf 2.1.} Soit $X$ un v.a. dans $\rd$ ayant la queue \`a
variation r\'eguli\`ere et $(\|X\|, \frac{X}{\|X\|})$ sa d\'ecomposition polaire. On va d'abord s'int\'eresser aux transformations qui ne changent que la partie sph\'erique de $X$. Plus exactement on prend une application mesurable $f:\sd \rightarrow \sd$ et on d\'efinit le nouveau vecteur $Y=\|X\|f(\frac{X}{\|X\|})$. Il est clair qu'en coordon\'ees polaires $Y=(\|X\|, f(\frac{X}{\|X\|}))$. Sous quelles conditions $Y$ reste-t-il encore le vecteur avec la queue r\'eguli\`ere? Quelle est la mesure spectrale de ce nouveau vecteur?
Le th\'eor\`eme suivant r\'epond \`a ces questions.

\begin{theoreme}\label{trans1}
Soit $X$ un v.a. dans $\rd$ tel que la condition de
variation r\'eguli\`ere (\ref{regulier}) a lieu avec l'exposant $\alpha$ et la mesure
spectrale $\sigma$. Soit $f$ une application $\sigma$-p.p. continue
sur $\sd$ \`a valeurs dans $\sd$, et $\mu$ est la mesure image d\'efinie par $\mu=\sigma f^{-1}$. Alors le v.a. transform\'e $Y=(\|X\|, f(\frac{X}{\|X\|}))$ a la
queue \`a variation r\'eguli\`ere de m\^eme exposant que $X$ et de la mesure spectrale $\mu$.
\end{theoreme}

\proof Prenons $B\in\mathcal{B}(\sd)$ tel que $\mu(\partial B)=0$. Alors
$\forall r>0$,
\begin{eqnarray*}
\frac{x^\alpha}{L(x)}\reip\left\{\frac{Y}{\|Y\|}\in
B,\|Y\|>x\right\}&=&\frac{x^\alpha}{L(x)}\reip\left\{f\left(\frac{X}{\|X\|}\right)\in B,\|X\|>x\right\}\\
&=&\frac{x^\alpha}{L(x)}\reip\left\{\frac{X}{\|X\|}\in f^{-1}(B),\|X\|>x\right\}.
\end{eqnarray*}
Pour que le dernier terme converge vers $\mu(B)r^{-\alpha}$, il
suffit d'assurer gr\^ace \`a (\ref{def1}) que $\sigma(\partial f^{-1}(B))=0$. Notons $D$
l'ensemble des discontinuit\'es de $f$, alors $\sigma(D)=0$. On va
montrer que \begin{equation}\label{baohan}
\partial f^{-1}(B)\subset f^{-1}(\partial B)\cup D.
\end{equation}
Si $x\in\partial f^{-1}(B)\backslash D$ alors $f$ est continue en
$x$. Comme $x\in\partial f^{-1}(B)$, il existe deux suites $x_n\in
f^{-1}(B)$ et $y_n\in f^{-1}(B)^\complement$, $n=1,2,\ldots$, telles
que $x_n\rightarrow x$ et $y_n\rightarrow x$. On en d\'eduit
\[f(x_n)\rightarrow f(x), \  f(x_n)\in B\]
et
\[f(y_n)\rightarrow f(x), \  f(y_n)\notin B.\]
Cela implique $f(x)\in\partial B$, ainsi $x\in f^{-1}(\partial B)$ d'o\`u (\ref{baohan}).
Puisque $\mu(\partial B)=\sigma f^{-1}(\partial B)=0$, on en d\'eduit
\[\sigma(\partial f^{-1}(B))\leq\sigma f^{-1}(\partial B)+\sigma(D)=0\]
qui compl\`ete la d\'emonstration. \qed

\begin{corollaire}\label{cortrans1}
Soit $X$ un v.a. dans $\Bbb{R}^2$ qui v\'erifie la
condition de variation r\'eguli\`ere (\ref{regulier}) avec l'exposant $\alpha$ et la mesure spectrale $\sigma$. On identifie $S^{1}$ avec l'intervalle $ [0,2\pi)$ et on suppose que la mesure spectrale de $X$ est uniforme, c'est-\`a-dire $d\sigma /d \theta =1/2\pi, \   \theta\in [0,2\pi)$. Soit $\mu$ une mesure de probabilit\'e sur $\sd$ avec la fonction de r\'epartition $F(x)=\mu ([0,x]),x\in [0,2\pi)$. Si $Y$ est un
v.a. d\'efini par $Y=(\|X\|, F^{-1}(\frac{X}{2\pi\|X\|}))$ o\`u $F^{-1}$ est la fonction de quantile correspondant \`a $F$, alors $Y$ v\'erifie la
condition de r\'egularit\'e avec la mesure spectrale $\mu$.
\end{corollaire}

\par\rm{\bf Remarque 1.}\quad Ce corollaire montre avec \'evidence l'utilit\'e des applications du th\'eor\`eme \ref{trans1} au probl\`eme de simulation.
\vspace{0.5cm}

\par\rm{\bf Remarque 2.}\quad La condition de continuit\'e de $f$ $\sigma$-p.p. est importante. L'exemple suivant montre que le r\'esultat du th\'eor\`eme \ref{trans1} n'est plus vrai si l'on omet cette condition.
\vspace{0.5cm}

\par\rm{\bf Exemple 1.}\quad Soient $F_{1}$ et $F_{2}$ deux lois d\'efinies sur $\mathbb{R}_+$ telles que 
\begin{enumerate}
\item la loi $F_{i}$ n'a pas de queue \`a variation r\'eguli\`ere, $i=1,2$,

\item la loi $F=\frac{1}{2}(F_{1}+F_{2})$ a la queue \`a variation r\'eguli\`ere avec l'exposant $\alpha$ et les constantes de normalisation $b_{n}$.
\end{enumerate}
Soient $\{1/n\}$ et $\{-1/n\}$ deux suites des points sur la sph\`ere unit\'e $S^1=(-\pi,\pi]$, not\'ees $\{x_{n}^+\}$ et $\{x_{n}^-\}$. Soit $l_{n}^\pm$ la demi-droite sortant de $0$ et passant par le point $x_{n}^\pm$. D\'efinissons deux suites des segments $\{\Delta_{n}^+\}$ et $\{\Delta_{n}^-\}$ par $\Delta_{n}^\pm=l_{n}^\pm\cap((-\pi,\pi]\times[n,n+1))$. Soit $\reip$ la loi sur $(-\pi,\pi]\times(\mathbb{R}_+\backslash\{0\})$ d\'efinie de la fa\c con suivante: son support est la r\'eunion de tous les intervalles $\Delta_{i}^\pm, \  i=1,2,\ldots$, et la restriction de $\reip$ sur $\bigcup_{i=1}^\infty \Delta_{i}^+$ (respectivement sur $\bigcup_{i=1}^\infty \Delta_{i}^-$) co\"{\i}ncide avec $1/2F_{1}$ transf\'er\'ee sur $l_{n}^+$ (respectivement avec $1/2F_{2}$ transf\'er\'ee sur $l_{n}^-$). Soit $X$ un v.a. de la loi $\reip$. On v\'erifie facilement que $\reip$ est la loi ayant la queue r\'eguli\`ere avec $b_{n}$ comme les constantes de normalisation dont la mesure spectrale est $\sigma=\delta_{\{0\}}$. 

Si l'on d\'efinit $f: (-\pi,\pi]\rightarrow(-\pi,\pi]$ par
\[f(x)=\left\{\begin{array}{ll}\frac{\pi}{2}&x>0,\\0&x=0,\\ -\frac{\pi}{2}&x<0,\end{array}\right.\]
on obtient que $\sigma f^{-1}=\sigma$ tandis que $Y=(\|X\|,f(\frac{X}{\|X\|}))$ n'a pas de queue r\'eguli\`ere. R\'eellement, par exemple, on a pour $\varepsilon>0$
\begin{eqnarray*}
n\reip\left\{\frac{Y}{\|Y\|}\in(\frac{\pi}{2}-\varepsilon,\frac{\pi}{2}+\varepsilon),\|Y\|>rb_{n}\right\}&=&n\reip\left\{\frac{X}{\|X\|}\in(0,\pi],\|X\|>rb_{n}\right\}\\&=&n(1-F_{1}(rb_{n}))
\end{eqnarray*}
qui ne converge nulle part par le choix de $F_{1}$.
\vspace{0.5cm}

{\bf 2.2.}  On consid\`ere maintenant les transformations ne modifiant que la partie radiale du vecteur initial. Etant donn\'e une fonction $h: \sd\rightarrow\mathbb{R_+}$ on d\'efinit le nouveau vecteur al\'eatoire $Y=Xh(\frac{X}{\|X\|})=(\|X\|h(\frac{X}{\|X\|}), \frac{X}{\|X\|})$. Le r\'esultat suivant donne les conditions sous lesquelles la propri\'et\'e de r\'egularit\'e de queue soit pr\'eserv\'ee.

\begin{theoreme}\label{trans2}
Soit $X$ un v.a. dans $\rd$ ayant la queue
r\'eguli\`ere d'exposant $\alpha$ avec la mesure spectrale $\sigma$. Soit $h$ une fonction
$\sigma$-p.p. continue et born\'ee sur $\sd$ \`a valeurs dans
$\Bbb{R}_+$, $\mu$ une mesure finie sur $\sd$ avec la densit\'e
$h(x)^\alpha$ par rapport \`a $\sigma$. Alors le v.a.
$Y=(\|X\|h(\frac{X}{\|X\|}),
\frac{X}{\|X\|})$ a la queue \`a variation
r\'eguli\`ere de m\^eme exposant $\alpha$ et de mesure spectrale $\mu$.
\end{theoreme}

La d\'emonstration est report\'ee dans la section suivante.

Les deux contre-exemples ci-dessous montrent que la condition que $h$ est $\sigma$-p.p. continue et born\'ee est r\'eellement importante pour pr\'eserver la r\'egularit\'e. L'exemple 2 pr\'esente une fonction $h$ $\sigma$-p.p. continue mais non-born\'ee pour laquelle le r\'esultat du th\'eor\`eme \ref{trans2} n'a pas lieu, tandis que la fonction $h$ de l'exemple 3 sera born\'ee mais non $\sigma$-p.p. continue.
\vspace{0.5cm}

\par\rm{\bf Exemple 2.}\quad {\it D\'efinition de $X$.}  Soit $\tau$ une mesure discr\`ete sur $S^1$ d\'efinie par
\[\tau(\{b_k\})=q_k=\frac{1}{k(k+1)}, \  b_k=\pi-\frac{\pi}{2^{k-1}}, \  k=1,2,\ldots .\]
Il est clair que $\sum\limits_{k=1}\limits^\infty q_k=1$ et
$b_k\in[0,\pi)$. Notons $L_k$ la demi-droite sortant de $0$ et passant par le point
$b_k$, c'est-\`a-dire $L_k=\{cb_k,c>0\}.$ Soit $Q_k$ une mesure sur $L_k$ dont la fonction de r\'epartition
est d\'efinie par
\[
F_k(x)=\left\{
\begin{array}{ll}
0&0<x<1\\
1-k^{-\nu}x^{-\alpha}&x\geq 1
\end{array}\right. \ 
\]
o\`u $\nu>0, \  k=1,2,\ldots$. Supposons que $X$ soit un v.a. dans $\Bbb{R}^2$
de la loi $\reip$ d\'efinie par:
\[\reip (A)=\sum_{k=1}^\infty q_kQ_{k}(A\cap L_k), \  A\in\mathcal{B}(\Bbb{R}^2).\]
 D\'efinissons la mesure
$\sigma$ sur $S^1$ par 
\begin{equation}\label{defsig} \sigma
(B)=\sum_{\{k|b_k\in B\}}q_kk^{-\nu}, \   B\in\mathcal{B}(S^1).
\end{equation} 
Cette mesure est bien d\'efinie car $\sum\limits_{k=1}^\infty q_kk^{-\nu}< 1$.
Maintenant pour tous $B\in\mathcal{B}(S^1)$ avec $\sigma (\partial B)=0$ et $\forall r>1$ on a
\begin{eqnarray}
r^\alpha\reip\left\{\frac{X}{\|X\|}\in B, \|X\|>r\right\}&=&\sum_{\{k|b_k\in
B\}}r^\alpha\reip\left\{\frac{X}{\|X\|}=b_k, \|X\|>r\right\}\nonumber\\
&=&r^\alpha\sum_{\{k|b_k\in B\}}q_k(1-F_k(r))\nonumber\\ &=&r^\alpha
r^{-\alpha}\sum_{\{k|b_k\in B\}}q_kk^{-\nu}\nonumber\\
&=&\sigma(B).\label{defsig2}
\end{eqnarray}
Cela signifie que $X$ a la loi avec la queue r\'eguli\`ere d'exposant $\alpha$ et de mesure spectrale $\sigma$.

On passe \`a la construction de notre function $h$. Prenons les intervalles $I_k,$ $k=1,2,\ldots$ sur
$S^1=[0,2\pi)$
\[I_1=(\frac{7\pi}{4}, 2\pi)\cup [0,\frac{\pi}{4}), \  I_k=(b_k-\frac{\pi}{2^{k+1}}, \   b_k+\frac{\pi}{2^{k+1}}), \   k\geq 2.\]
Puisque la distance entre $b_k$ et $b_{k+1}$ est
${\displaystyle\frac{\pi}{2^k}}$, les intervalles $I_1, I_2,\ldots$
sont disjoints et $b_k\in I_k$ pour chaque $k$. Notre fonction $h$
est d\'efinie par 
\begin{equation}\label{defh}
h(x)=\sum_{k=1}^\infty
k^\beta\one_{I_k}(x), \end{equation}
o\`u $\beta$ est tel que  $\frac{1}{\alpha}<\beta<\frac{1+\nu}{\alpha}$. \'Evidemment $h$ est
$\sigma$-p.p. continue et non-born\'ee. Si
$Y=Xh(\frac{X}{\|X\|})$, alors pour $\forall r>1$
\begin{eqnarray*}
r^\alpha\reip\left\{\frac{Y}{\|Y\|}\in S^1, \|Y\|>r\right\}&=&r^\alpha\reip\left\{\|X\|h\left(\frac{X}{\|X\|}\right)>r\right\}\\
&=&r^\alpha\sum_{k=1}^\infty\reip\left\{\frac{X}{\|X\|}=b_k,
\|X\|>\frac{r}{k^\beta}\right\}\\
&\geq&r^\alpha\sum_{\{k|\frac{r}{k^\beta}<
1\}}\frac{1}{k(k+1)}\\
&\geq&\frac{r^\alpha}{r^{1/\beta}+1},
\end{eqnarray*}
d'o\`u suit la convergence vers l'infini quand $r\rightarrow\infty$ du terme \`a droite, ce qui n'aurait pas lieu si le th\'eor\`eme \ref{trans2} \'etait applicable.
\vspace{0.5cm}

\par\rm{\bf Remarque 3.}\quad En vu du th\'eor\`eme \ref{cdmoinsforte} on pourrait penser que la condition suivante et moins restrictive
\[\exists \delta>0 \  \mbox{tel que} \  \int_{S^1}h^{\alpha+\delta}d\sigma<\infty\]
sera suffisante pour pr\'eserver la r\'egularit\'e de queue. Notre exemple montre que ce n'est pas le cas. R\'eellement, si $\delta$ est suffisamment petit,
\[\int_{S^1}h^{\alpha+\delta}(\theta)\sigma(\theta)=\sum_{k=1}^{\infty}k^{(\alpha+\delta)\beta-\nu}q_{k}\leq\sum_{k=1}^{\infty}k^{(\alpha+\delta)\beta-\nu-2}<\infty.\]

\par\rm{\bf Exemple 3.}\quad {\it D\'efinition de $X$.} On consid\`ere la fonction $g(x)$ sur $\Bbb{R}_+$
\[g(x)=\sum_{k=0}^\infty\frac{1}{2^k}\one_{(k,k+1]}(x).\]
Notons son graphe par $D_g$:
\[D_g=\{(x,y)| y=g(x), x>0\}.\]
Soit $Q$ une mesure sur $\Bbb{R}_+$ dont la fonction de
r\'epartition est d\'efinie par
\[F_Q(x)=1-G(x)=1-1 \wedge x^{-\alpha}.\]
L'application $\pi$ de $\Bbb{R}_+$ \`a $\mathbb{R}^2$ d\'efinie par
\[\pi : r\mapsto (r,g(r))\]
transforme $Q$ en mesure image $Q\pi^{-1}$ qui sera concentr\'ee sur $D_{g}$.

Supposons que $X$ est un v.a. dans $\Bbb{R}^2$ de la loi $\reip$ suivante:
\[\reip(A)=\frac{1}{2}Q\pi^{-1}(A\cap D_g)+\frac{1}{2}Q(A\cap E), \;  \;   A\in\mathcal{B}(\Bbb{R}^2),\]
o\`u $E=\{(x,0)|x>0\}$.
On v\'erifie que $X$ satisfait (\ref{regulier2}). Prenons d'abord $B=[a,2\pi)$,
$a\in(0,2\pi)$. Notons $k_a=\min\{k |\frac{1}{2^k}\leq
a\}$, alors
\[((k_a,\infty)\times B) \cap D_g=\emptyset.\]
Donc pour tout $r>k_a$ on a
\begin{equation}\label{ex3rel1}
r^\alpha\reip\left\{\frac{X}{\|X\|}\in B, \   \|X\|>r\right\}=0.
\end{equation}
Ensuite, si $B=[0,a)$, pour $\forall r>k_{a}$
\begin{eqnarray}
r^\alpha\reip\left\{(r,\infty)\times B\right\}&=&r^\alpha\reip\{(r,\infty)\times (0,a)) \cup((r,\infty)\times \{0\}\})\nonumber\\
&=& r^\alpha(\frac{1}{2}Q\pi^{-1}(((r,\infty)\times (0,a)) \cap
D_g)+\frac{1}{2}Q((r,\infty)))\nonumber\\
&=& r^\alpha G(r)=1.\label{ex3rel2}
\end{eqnarray}
Les relations (\ref{ex3rel1}), (\ref{ex3rel2}) donnent (\ref{regulier}) avec $\sigma=\delta_{\{0\}}$.

{\it D\'efinition de $h$.} On pose $h(x)=\one_{(0,2\pi)}(x)$. Alors l'ensemble des discontinuit\'es de $h$ sera $\{0\}$ et puisque $\sigma(\partial\{0\})=\sigma (\{0\})=1$, la fonction $h$ n'est pas $\sigma$-p.p. continue.

Par les arguments analogues aux pr\'ec\'edents on trouve que le vecteur $Y=Xh(\frac{X}{\|X\|})$ satisfait la condition de r\'egularit\'e avec la mesure spectrale $\mu=\frac{1}{2}\delta_{\{0\}}$. Par cons\'equent $d\mu/d\sigma\neq h$.
\vspace{0.5cm}

Si les variables al\'eatoires $\frac{X}{\|X\|}$ et $\|X\|$ sont ind\'ependantes il y a une
condition moins forte sur $h$ telle que la r\'egularit\'e soit pr\'eserv\'ee sous la transformation.

\begin{theoreme}\label{cdmoinsforte}
Soit $X$ un v.a. dans $\mathbb{R}^d$ satisfaisant la
condition (\ref{regulier}) avec l'exposant $\alpha$ et la mesure spectrale $\sigma$, $h$ une
fonction d\'efinie sur $S^{d-1}$ \`a valeurs dans $\mathbb{R}_+$
telle que $\int_{S^{d-1}}h^{\alpha+\varepsilon}d\sigma <\infty$ pour
un $\varepsilon >0$. Si les variables $
\frac{X}{\|X\|}$ et $\|X\|$ sont ind\'ependantes, alors le vecteur
transform\'e $Y=Xh(\frac{X}{\|X\|})=(\|X\|h(\frac{X}{\|X\|}), \  
 \frac{X}{\|X\|})$ v\'erifie la condition
(\ref{regulier}) avec la mesure spectrale $\mu$ telle que
$d\mu/d\sigma=h^\alpha$.
\end{theoreme}

\proof En repr\'esentant $Y$ sous la forme $Y=\|X\|Z$, o\`u $Z=\frac{X}{\|X\|}h(\frac{X}{\|X\|})$, on remarque que le r\'esultat suit directement du Th. 4.15 \cite{resmov}. En connection de ce r\'esultat on peut mentionner \cite{breiman} (pour $d=1$) et \cite{davydov2000flt}, Lemme 3.9.
\qed

En randomisant la fonction $h$, on d\'eduit imm\'ediatement du th\'eor\`eme \ref{cdmoinsforte} le corollaire suivant.
\begin{corollaire}
Supposons que $X$ satisfait la condition (\ref{regulier}) avec l'exposant $\alpha$ et la mesure spectrale $\sigma$.  Soit $\{Z(\theta),\theta\in\sd\}$ un processus stochastique ind\'ependant de $X$ dont les trajectoires sont presque s\^urement positives et $\sigma$-p.p. continues. Si pour un $\varepsilon>0$ 
\[\int_{\sd}\reie (Z^{\alpha+\varepsilon}(\theta))\sigma(d\theta)<\infty\]
 alors le vecteur $Y=XZ(\frac{X}{\|X\|})$ a la queue r\'eguli\`ere de m\^eme exposant $\alpha$ que $X$ et de mesure spectrale $\mu$ telle que $d\mu/d\sigma=\reie (Z(\theta))^\alpha$.
\end{corollaire}

\section{Preuves}

\par\rm{\bf D\'emonstration du th\'eor\`eme \ref{trans2}}\quad
Rappelons que $Q_n$, $Q$ sont les mesures associ\'ees avec $X$ et d\'efinies par (\ref{defQn})
et (\ref{defQ}). Soit $\tilde{Q}_{n}$, $\tilde{Q}$ les mesures associ\'ees avec le vecteur $Y=Xh(\frac{X}{\|X\|})$ et d\'efinies de la m\^eme fa\c{c}on, i.e. 
\begin{equation}\label{deftildeQnQ}
\tilde{Q}_n((r,\infty)\times B)=n\reip\left\{\frac{Y}{\|Y\|}\in B, \|Y\|>rb_n\right\},  \; \;  \tilde{Q}=m_\alpha\times\mu,
\end{equation}
o\`u $B\in\mathcal{B}(\sd)$ et $r>0$, $b_{n}$ est le m\^eme que dans (\ref{defQn}). Notons $\mathcal{S}_{+}$ la famille des fonctions sur $\rd$ avec les supports qui sont s\'epar\'es du point z\'ero, i.e.
\[\mathcal{S}_{+}=\{f | \ \exists \varepsilon>0 \ \mbox{tel que} \  supp(f)\subset (\varepsilon,\infty)\times\sd\}.\]
Pour d\'emontrer le th\'eor\`eme, d'apr\`es le th\'eor\`eme de portmanteau, il nous suffit d'\'etablir la convergence 
\begin{equation}\label{cvtildeQ}
\int fd\tilde{Q}_n\rightarrow\int fd\tilde{Q},  \;  \;  n\rightarrow\infty,
\end{equation}
pour toute $f\in \mathcal{S}_{+}$ continue et born\'ee. D\'efinissons l'application $\varphi$ par
\begin{eqnarray*}
\varphi:\mathbb{R}_+\times\sd & \rightarrow&\mathbb{R}_+\times\sd,\\
(\rho,\theta)&\mapsto & (\rho h(\theta),\theta).
\end{eqnarray*}
On remarque que $Y=\varphi(X)$, et 
\[Q_n\varphi^{-1}((r,\infty)\times B)=n\reip\left\{\frac{Y}{\|Y\|}\in B, \|Y\|>rb_n\right\},\]
 \begin{eqnarray*} Q\varphi^{-1}((r,\infty)\times B)&=&Q\{(\rho,\theta)|\theta\in B,\rho h(\theta)\in(r,\infty)\}\\
&=& \int_B\int\one_{(r,\infty)} (\rho
h(\theta))m_\alpha (d\rho)\sigma(d\theta)\\
&=&\int_B\sigma(d\theta)\alpha\int\one_{(r,\infty)}(\rho
h(\theta))\rho^{-\alpha-1}d\rho\\
&=&\mu(B)r^{-\alpha}.\end{eqnarray*}
Par cons\'equent,
\[\int fd\tilde{Q}_n=\int(f\circ\varphi)dQ_{n},\]
\[\int fd\tilde{Q}=\int(f\circ\varphi)dQ.\]
La fonction $f\circ\varphi$ est born\'ee et $Q$-p.p. continue gr\^ace \`a $Q$-p.p. continuit\'e de $\varphi$. De plus, $f\circ\varphi\in \mathcal{S}_{+}$ puisque $h$ est suppos\'ee born\'ee. En effet si $h\leq M$ et $\mbox{supp}(f)\subset (\varepsilon,\infty)\times\sd$, alors $\mbox{supp}(f\circ\varphi)\subset (\varepsilon/M,\infty)\times\sd$.
Car $X$ a la queue r\'eguli\`ere de l'exposant $\alpha$ avec la mesure spectrale $\sigma$, gr\^ace \`a la convergence \'equivalente (\ref{def4}), on a pour toute $f\in  \mathcal{S}_{+}$
\[ \int(f\circ\varphi)dQ_{n}\rightarrow \int(f\circ\varphi)dQ\]
ce qui donne (\ref{cvtildeQ}).
  \qed

\bibliographystyle{plain}

\end{document}